\documentclass{article}
\usepackage{amsmath,amsthm}
\usepackage{graphicx,epstopdf,epsfig,multirow,epic,bm}

\normalsize \rm
\begin{document}

\begin{center}
{\Large  \textbf {Exact solutions for geodesic distance on treelike models with some constraints}}\\[12pt]
{\large \quad Xudong Luo$^{a,}$\footnote{~The author's E-mail: luoxudong117@163.com.},\quad  Fei Ma$^{b,}$\footnote{~The corresponding author's E-mail: mafei123987@163.com. } \; and \; Wentao Xu$^{c,}$\footnote{~The author's E-mail: xuwentao1216@163.com. } }\\[6pt]
{\footnotesize $^{a}$ College of Mathematics and Statistics, Northwest Normal University, Lanzhou  730070, China\\
 $^{b}$ School of Electronics Engineering and Computer Science, Peking University, Beijing 100871, China\\
$^{c}$ School of Earth Sciences and Engineering, Nanjing University, Nanjing 210023, China}\\[12pt]
\end{center}

\begin{quote}
\textbf{Abstract:}  Geodesic distance, commonly called shortest path length, has proved useful in a great variety of disciplines. It has been playing a significant role in search engine at present and so attracted considerable attention at the last few decades, particularly, almost all data structures and corresponding algorithms suitable to searching information generated based on treelike models. Hence, we, in this paper, study in detail geodesic distance on some treelike models which can be generated by three different types of operations, including first-order subdivision, ($1,m$)-star-fractal operation and $m$-vertex-operation. Compared to the most best used approaches for calculating geodesic distance on graphs, for instance, enumeration method and matrix multiplication, we take useful advantage of a novel method consisting in spirit of the concept of vertex cover in the language of graph theory and mapping. For each kind of treelike model addressed here, we certainly obtain an exact solution for its geodesic distance using our method. With the help of computer simulations, we confirm that the analytical results are in perfect agreement with simulations. In addition, we also report some intriguing structure properties on treelike models of two types among them. The one obeys exponential degree distribution seen in many complex networks, by contrast, the other possesses all but leaf vertices with identical degree and shows more homogeneous topological structure than the former. Besides that, the both have, in some sense, self-similar feature but instead the latter exhibits fractal property. 

\end{quote}

\vskip 1cm

\section{INTRODUCTION}

The last several decades have witnessed an upsurge in complex network study which ranges from statistical physics, discrete applied mathematics, theoretical computer science, to biological science and chemistry, even to social science, and so forth. As a newborn mathematical tool, complex network has been proved useful in understanding most complex systems in both nature and real-life world. Included example networks have World Wide Web (WWW), citation networks, metabolic networks, protein-protein interaction networks and predator-prey webs \cite{Michael-2019}-\cite{Ulrich-2019}. In particular, one type of networked models have attracted more attention according to their own specific topological structure like tree. One of such examples is citation networks of scientific papers \cite{Julia-2019}. In fact, there have been a great deal of treelike models proposed to model many complex networks. Hence a number of their topological structure properties have been reported, for instance, geodesic distance \cite{Deng-2019}, mean first-passage time for random walk \cite{Peng-2018}, fractal phenomena \cite{Nobutoshi-2019} and so on.

Generally speaking, there are two mainstreams in the studies of complex networks at present. The one is to devote to probing the generation mechanisms driving the evolution of complex networks over time. The best known is the preferential attachment mechanism \cite{Albert-1999}. The other is to understand dynamic and function taking place on complex networks themselves. The former outlines ``why" and however the latter focuses on ``how". In this paper, our discussion falls into the region of the latter. A crucial issue included in the latter is to reveal and understand how the underlying topological structure of complex networks influences those dynamical behaviors over them. To do this, recent great efforts have been made. Among of which, geodesic distance is an important topological index of complex networks. Geodesic distance, commonly called shortest path length, is defined as the length of path between a pair of vertices. It is, in essence, not a fresh measure and has been widely studied in many disciplines, such as, graph theory. Despite that, the considerable importance and usefulness of geodesic distance have been highlighted by a variety of practical applications, such as, searching information on internet \cite{Spyros-2017}, signal integrity in communication networks, disease spreading on relationship networks among individuals \cite{Jia-2018}, navigation in spatial networks \cite{Wei-2014}, to name but a few.

As mentioned above, a central problem to address is how to analytically capture solutions for computation of geodesic distance on treelike models according to theoretical value and practical applications. Some relevant researches have been shown in some fields \cite{Zhang-2010,Zhang-2011}. Therefore, exact solutions for several types of treelike models have been obtained by taking advantage of some mature methods in which the best used is Laplacian spectra and eigenvectors of underlying structure from spectral graph theory. Along the line of such researches, a lot of results have been reported. Nonetheless, we here do not employ methods of this type and instead introduce some novel methods based on topological structure of treelike models themselves sufficiently. Although the nature of both our methods and that of \cite{Zhang-2010} is to built up a group of equations in an iterative manner, ours are, in some cases, more convenient to manipulate than the latter, at least on treelike models proposed in this paper.

This paper can be organized by the next several Sections. In Section 2, we will introduce some helpful definitions, including first-order subdivision, ($1,m$)-star-fractal operation and $m$-vertex-operation, in order to generate our research objects which are some well studied treelike models, and conventional notations, for example, vertex cover, surjection and bijection, to smoothly develop our main results. And then, three theorems as our main results, which are all based on a novel method for calculating geodesic distance on treelike models, are built up in Section 3. To show the practicality of our methods addressed here, in Section 4, we make use of them to derive exact solutions for geodesic distance on two classes of treelike models which coincide perfectly with that published results \cite{Zhang-2010,Zhang-2011}. Meanwhile, this convenience of our methods is well highlighted through reducing tough calculations in comparison with already established techniques, such as Laplacian spectra and eigenvectors. In addition, we also make a description of some interesting topological properties owned by the both treelike models themselves briefly. In conclusion, we close this paper by outlining the importance of our works and reporting some potential applications of the light shed by our methods to direct our future work.

\section{Definitions and notations}

It is conventional in graph theory terms to let $\mathcal{G}(\mathcal{V},\mathcal{E})$ be a graph whose vertex set and edge set are $\mathcal{V}$ and $\mathcal{E}$, respectively, and vertex number (order) and edge number (size) are denoted by $|\mathcal{V}|$ and $|\mathcal{E}|$ where symbol $||$ represents the cardinality of a set. Meanwhile, the notation $[1,n]$ is an integer set which consists precisely of those integers no more than $n$ and no less than $1$. For details to see Ref.\cite{Bondy-2008}.

\textbf{Definition 1} Given an arbitrary graph $\mathcal{G}(\mathcal{V},\mathcal{E})$, if one inserts a new vertex to every edge $uv\in \mathcal{E}$ then the resulting graph, denoted by $\mathcal{G'}_{1}(\mathcal{V'}_{1},\mathcal{E'}_{1})$, is called a subdivision of original graph $\mathcal{G}(\mathcal{V},\mathcal{E})$. Equivalently, such a subdivision can be obtained from graph $\mathcal{G}(\mathcal{V},\mathcal{E})$ by replacing every edge $uv\in \mathcal{E}$ by a unique path $uwv$ with length $2$ where internal vertex $w$ is that added vertex. Hereafter we think of such an operation as \textbf{\emph{subdivision}} and that resulting graph as \textbf{\emph{subdivision graph }}of graph $\mathcal{G}(\mathcal{V},\mathcal{E})$. More generally, such a subdivision is commonly regarded as \textbf{\emph{first-order subdivision }}vividly. It is worth noting that we in this paper focus mainly on how the first-order subdivision does make a considerable effect on geodesic distance of tree $\mathcal{T}(\mathcal{V},\mathcal{E})$. Here, Fig.1(b) shows a subdivision $\mathcal{T'}$ of tree $\mathcal{T}$ on seven vertices plotted in Fig.1(a).

For our purpose, it can immediately see using Def.1 that the subdivision graph $\mathcal{G'}_{1}(\mathcal{V'}_{1},\mathcal{E'}_{1})$ holds on $|\mathcal{V'}_{1}|=|\mathcal{V}|+|\mathcal{E}|$ and $|\mathcal{E'}_{1}|=2|\mathcal{E}|$. After $t$ time steps, the order $|\mathcal{V'}_{t}|$ and size $|\mathcal{E'}_{t}|$ of the subdivision graph $\mathcal{G'}_{t}(\mathcal{V'}_{t},\mathcal{E'}_{t})$ will obey a couple of equations as follows

\begin{equation}\label{Section-2-0}
|\mathcal{V'}_{t}|=|\mathcal{V}|+(2^{t}-1)|\mathcal{E}|, \qquad |\mathcal{E'}_{t}|=2^{t}|\mathcal{E}|.
\end{equation}

\textbf{Definition 2} Given an arbitrary graph $\mathcal{G}(\mathcal{V},\mathcal{E})$, if one not only inserts a vertex to every edge $uv\in \mathcal{E}$ but also connects $m$ other new vertices to this newly inserted vertex, then the resulting graph, denoted by $\mathcal{G^{\star}}_{1}(\mathcal{V^{\star}}_{1},\mathcal{E^{\star}}_{1})$, is called a\textbf{ \emph{(1,m)-star-fractal graph} }of original graph $\mathcal{G}(\mathcal{V},\mathcal{E})$. Specifically speaking, such a operation can be achieved from graph $\mathcal{G}(\mathcal{V},\mathcal{E})$ by directly inserting a star with $m$ leaves to every $uv\in \mathcal{E}$ and hence called\textbf{ \emph{(1,m)-star-fractal}}. It is obvious to say that the well known \textbf{\emph{T-fractal}} can be induced as a result of our ($1,m$)-star-fractal graph when parameter $m$ is supposed equal to $1$ \cite{Peng-2018}. As before, our aim is to study geodesic distance on ($1,m$)-star-fractal tree of tree $\mathcal{T}(\mathcal{V},\mathcal{E})$. An example as illustration of ($1,m$)-star-fractal of tree $\mathcal{T}$ with seven vertices is shown in Fig.1(c).

For brevity, based on Def.2, one can find out that ($1,m$)-star-fractal graph $\mathcal{G^{\star}}_{1}(\mathcal{V^{\star}}_{1},\mathcal{E^{\star}}_{1})$ has $|\mathcal{V^{\star}}_{1}|=|\mathcal{V}|+(1+m)|\mathcal{E}|$ vertices and $|\mathcal{E^{\star}}_{1}|=(2+m)|\mathcal{E}|$ edges. Similarly, for time step $t$, both vertex number and edge number of ($1,m$)-star-fractal graph $\mathcal{G^{\star}}_{t}(\mathcal{V^{\star}}_{t},\mathcal{E^{\star}}_{t})$, respectively, obey

\begin{equation}\label{Section-2-1}
|\mathcal{V^{\star}}_{t}|=|\mathcal{V}|+((2+m)^{t}-1)|\mathcal{E}|, \qquad |\mathcal{E^{\star}}_{t}|=(2+m)^{t}|\mathcal{E}|.
\end{equation}

\textbf{Definition 3} Given an arbitrary graph $\mathcal{G}(\mathcal{V},\mathcal{E})$, if one only connects $m$ new vertices to every vertex $u\in \mathcal{V}$ then the resulting graph, denoted by $\mathcal{G^{\odot}}_{1}(\mathcal{V^{\odot}}_{1},\mathcal{E^{\odot}}_{1})$, is called a\textbf{ \emph{m-vertex-operation graph} }of original graph $\mathcal{G}(\mathcal{V},\mathcal{E})$. Different from the above manipulations introduced in Def.1 and Def.2, operation here is manipulated on vertex and thus viewed as\textbf{ \emph{m-vertex-operation}} which leads to $|\mathcal{V}|$ stars where each includes one vertex $u\in \mathcal{V}$ as its central vertex and other $m$ new vertices attached to vertex $u$. At the same time, our goal is to put insight into calculating analytically geodesic distance on $m$-vertex-operation tree of tree $\mathcal{T}(\mathcal{V},\mathcal{E})$. Such an example is plotted in Fig.1(d) serving as an illustration of $m$-vertex-operation of tree $\mathcal{T}$ with seven vertices.

For convenience, one can understand on the basis of Def.3 that in $m$-vertex-operation graph $\mathcal{G^{\odot}}_{1}(\mathcal{V^{\odot}}_{1},\mathcal{E^{\odot}}_{1})$ $|\mathcal{V^{\odot}}_{1}|$ is equal to $(m+1)|\mathcal{V}|$ and $|\mathcal{E^{\odot}}_{1}|$ equal to $m|\mathcal{V}|+\mathcal{E}$. After $t$ time steps, the  $m$-vertex-operation graph $\mathcal{G^{\odot}}_{t}(\mathcal{V^{\odot}}_{t},\mathcal{E^{\odot}}_{t})$ holds on

\begin{equation}\label{Section-2-2}
|\mathcal{V^{\odot}}_{t}|=(m+1)^{t}|\mathcal{V}|, \qquad |\mathcal{E^{\odot}}_{t}|=((m+1)^{t}-1)|\mathcal{V}|+|\mathcal{E}|.
\end{equation}

Note also that for a given graph $\mathcal{G}(\mathcal{V},\mathcal{E})$ the three resulting graphs mentioned above all have the same average degree $\langle k\rangle\sim 2$, which is typically written as $\langle k\rangle=2|\mathcal{E}|/|\mathcal{V}|$, in the limit of large graph size. It is straightforward to say that graphs with an identical average degree may have completely different topological structure from one another. Suppose that original seed graph $\mathcal{G}(\mathcal{V},\mathcal{E})$ is, for instance, a connected graph with no leaves then the subdivision graph $\mathcal{G'}_{t}(\mathcal{V'}_{t},\mathcal{E'}_{t})$ still remains connected and has no leaves but the latter two resulting graphs, $\mathcal{G^{\star}}_{t}(\mathcal{V^{\star}}_{t},\mathcal{E^{\star}}_{t})$ and $\mathcal{G^{\odot}}_{t}(\mathcal{V^{\odot}}_{t},\mathcal{E^{\odot}}_{t})$, will have various types of treelike branches. Besides that, the latter both are very likely to follow different properties from each other in some other respects as we will show shortly.

\begin{figure}
\centering
  \includegraphics[height=7cm]{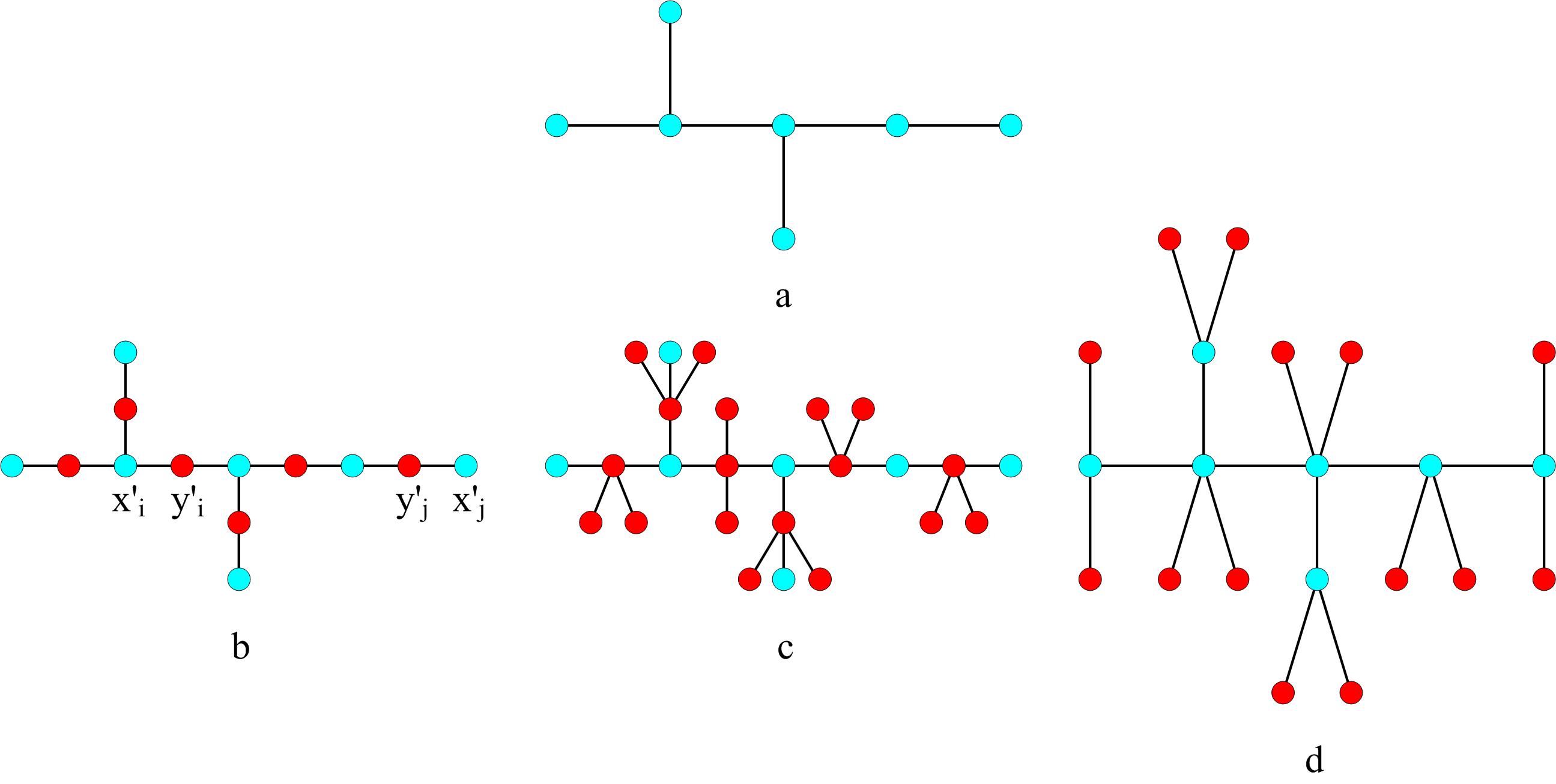}\\
{\small Fig.1. The diagrams of a tree $\mathcal{T}$ on seven vertices and the resulting trees obtained from tree $\mathcal{T}$ by applying three different types of operations. Panel (a) shows tree $\mathcal{T}$. By applying \emph{first-order-subdivision} to each edge of tree $\mathcal{T}$ results in a new tree as plotted in panel (b). Similarly, panel (c) shows a resulting tree using \emph{(1,m)-star-fractal} on each edge of tree $\mathcal{T}$ where $m=2$. The final panel (i.e., panel (d)) describes a tree generated by implementing \emph{m-vertex-operation} on each vertex of tree $\mathcal{T}$ where $m=2$.        }
\end{figure}

Below will introduce a helpful definition used later to accomplish our proofs smoothly that is in practice a product stem from vertex cover in the jargon of graph theory. For detail to see \cite{Bondy-2008}, we here just discuss such problems in the simplest and most fundamental situation.

\textbf{Definition 4} Given a path $\mathcal{P}$ with $n$ vertices, it is easy to see that the $n$ vertices can be grouped into two disconnected vertex sets, without loss of generality, written as set $X=\{x_{1},x_{2},...,x_{\lceil n/2\rceil}\}$ and set $Y=\{y_{1},y_{2},...,y_{\lfloor n/2\rfloor}\}$. The above classification of vertices is in fact a bipartition of vertex set. More generally, vertices of set $X$ and vertices of set $Y$ can be alternatively arranged on the path $\mathcal{P}$ in a reasonable manner such that arbitrary vertex pair $x_{i}$ and $x_{j}$ is not connected directly by an edge and similarly for all pairs of vertices $y_{i}$ and $y_{j}$. This suggests that either vertex set $X$ or vertex set $Y$ is a \textbf{\emph{vertex cover}} of path $\mathcal{P}$ where $|X|=|Y|+1$ when $n$ is odd and $|X|=|Y|$ otherwise. The both vertex sets will alternatively paly a crucial role on the growth process of building up our main results in the rest of this paper.

Now let us take a common yet important terminology from analysis mathematics which has been utilized in a great variety of science fields.

\textbf{Definition 5} Given two sets $X$ and $Y$ that all consist of at least one element, one can introduce a \textbf{\emph{mapping}} $f$ from $X$ to $Y$ such that for a provided element $x$ of set $X$ there must be a unique element $y$ belonging to set $Y$ meeting $f(x)=y$. This can be simply expressed in the following mathematically

$$\forall x\in X, \quad \exists y\in Y, \quad s.t.,\quad f:\; x\mapsto \; y$$
where the image set of set $X$ may be written as $f_{X\mapsto Y}=\{y|f(x)=y, \; y\in Y\}$. We is in this paper interested in the below two kinds of mappings between sets $X$ and $Y$.

\emph{case 1} If both $|X|>|Y|$ and $|f_{X\mapsto Y}|=|Y|$ are true, then this mapping $f$ is considered \textbf{\emph{surjection}}. Besides, for each element $y$ of image set $Y$, there exist $n$ distinct pre-images $x_{i}$ ($i\in [1,n]$), i.e., $f^{-1}(y)=x_{i}$, then the surjection $f$ is considered \textbf{\emph{n-regular}}. It is clear to the eye that both sets $X$ and $Y$ satisfy $|X|=n |Y|$ in question.

\emph{case 2} If the surjection $f$ under consideration holds $|X|=|Y|$ then it can be thought of as a \textbf{\emph{bijection}}, also called one-one mapping.

Meanwhile, the compound mapping between mappings $f$ and $g$ can be expressed as $f\circ g$ mathematically.

By far, we have introduced some helpful definitions and notations used later. As stated above, the topics of this paper focus principally on many discussions correlated to geodesic distance on treelike models with some constraints. Now, let us turn our insight to these problems.

\section{Main results}

We will in this section show our main results along with corresponding brief proofs which are organized by three theorems in form.

\textbf{Theorem 1} Given an arbitrary tree $\mathcal{T}(\mathcal{V},\mathcal{E})$, the exact solution for geodesic distance $\mathcal{S'}$ of its first-order subdivision tree $\mathcal{T'}(\mathcal{V'},\mathcal{E'})$ is

\begin{equation}\label{Section-3-1-0}
\mathcal{S'}=8\mathcal{S}-2|\mathcal{V}|(|\mathcal{V}|-1)
\end{equation}
in which $\mathcal{S}$ is a known expression to geodesic distance of tree $\mathcal{T}(\mathcal{V},\mathcal{E})$.

\textbf{\emph{Proof}} At first, suppose that the geodesic distance on tree $\mathcal{T}(\mathcal{V},\mathcal{E})$ is $\mathcal{S}$. After applying first-order subdivision to each edge of tree $\mathcal{T}(\mathcal{V},\mathcal{E})$, the first-order subdivision tree $T'(\mathcal{V'},\mathcal{E'})$ will have two types of vertices, without loss of generality, which are divided into two disjoint vertex sets $X'$ and $Y'$. Set $X'$ is in essence set $\mathcal{V}$ and set $Y'$ is made up all the newly added vertices, namely, $Y'=\mathcal{V'}-\mathcal{V}$. In order to precisely calculate geodesic distance $\mathcal{S'}$, it is straightforward to compute three classes of geodesic distances, one for vertex pairs $<x'_{i},x'_{j}>$ of set $X'$, one for vertex pairs $<y'_{i},y'_{j}>$ of set $Y'$ as well the latter for vertex pairs $<x'_{i},y'_{j}>$ between set $X'$ and set $Y'$.

\emph{Case 1.1} For a given vertex pair $<x'_{i},x'_{j}>$ of set $X'$, there must be a bijection $f_{1}$ between set $X'$ and set $\mathcal{V}$ such that $f_{1}(<x'_{i},x'_{j}>)=<x_{i},x_{j}>$ here vertices $x_{i}$ and $x_{j}$ are in set $\mathcal{V}$. In fact, such a bijection $f_{1}$ is self-mapping and hence one can write

\begin{equation}\label{Section-3-1-1}
\mathcal{S'}(1)=2\mathcal{S}
\end{equation}
where $\mathcal{S'}(1)$ is the sum of distances of all possible vertex pairs in set $X'$.

\emph{Case 1.2} Different from case 1.1, there indeed is a self-mapping between set $Y'$ and set $\mathcal{V'}-\mathcal{V}$ but is no use for our problem. Taking into consideration results in case 1.1 known to us, the current issue is to build connection to Eq.(\ref{Section-3-1-1}). Therefore, for a given vertex pair $<y'_{i},y'_{j}>$ of set $Y'$, we should find out a bijection $f_{2}$ between set $Y'$ and set $X'$. It does not require much effort to do this according to statement in Def.4. To see why this is, let us pay attention on the both disjoint vertex sets $X'$ and $Y'$. For arbitrary vertex pair $<y'_{i},y'_{j}>$ of set $Y'$, there must be a unique path $\mathcal{P}_{y'_{i}y'_{j}}$ connecting vertices $y'_{i}$ and $y'_{j}$ whose vertices are alternatively chosen from the both sets above such that two vertices arranged at the two endpoints of path $\mathcal{P}_{y'_{i}y'_{j}}$ are $y'_{i}$ and $y'_{j}$. One immediately obtain an extension $\mathcal{P}_{x'_{i}x'_{j}}$ of this path $\mathcal{P}_{y'_{i}y'_{j}}$ by adding two edges $x'_{i}y'_{i}$ and $x'_{j}y'_{j}$, meaning which there is a mapping $f^{*}$ between vertex pairs $<y'_{i},y'_{j}>$ and $<x'_{i},x'_{j}>$, see Fig.1(b). Meanwhile, it is not hard to prove mapping $f^{*}$ to be bijection in terms of combinations among bijection $f_{1}$, first-order subdivision and tree itself. Thus, the compound function between the candidate mapping $f^{*}$ and bijection $f_{1}$ is chosen as our desired bijection $f_{2}$, i.e., $f_{2}=f_{1}\circ f^{*}$. Armed with these demonstrations, we have

\begin{equation}\label{Section-3-1-2}
\mathcal{S'}(2)=\mathcal{S'}(1)-2\frac{|\mathcal{V}|(|\mathcal{V}|-1)}{2}
\end{equation}
in which $\mathcal{S'}(2)$ is geodesic distance of all possible vertex pairs in set $Y'$.

\emph{Case 1.3} The remainder of our problem is to capture the expression of geodesic distance for all possible vertex pairs $<x'_{i},y'_{j}>$ between set $X'$ and set $Y'$. Along the research line of case 1.2, for a pair of vertices $x'_{i}$ and $y'_{j}$, there is a unique path $\mathcal{P}_{x'_{i}y'_{j}}$ which can also be expanded to another path $\mathcal{P}_{x'_{i}x'_{j}}$ by adding an additional edge $y'_{j}x'_{j}$ under a similar mapping $f^{**}$ to mapping $f^{*}$ in case 1.2. Therefore, we may create a mapping $f_{3}=f_{1}\circ f^{**}$ which must be a surjection between vertex pairs $<x'_{i},y'_{j}>$ and $<x'_{i},x'_{j}>$ but not bijection. One of reasons for this is another path $\mathcal{P}_{y'_{i}x'_{j}}$ may be reduced as the path $\mathcal{P}_{x'_{i}x'_{j}}$ by adding edge $x'_{i}y'_{i}$ as well. In a word, there are two distinct pre-images under surjection $f_{3}$, i.e., $<x'_{i},y'_{j}>$ and $<y'_{i},x'_{j}>$, such that $f^{-1}_{3}(<x'_{i},x'_{j}>)=<x'_{i},y'_{j}>=<y'_{i},x'_{j}>$. By Def.5, such a surjection is in principle 2-regular and so the geodesic distance $\mathcal{S'}(3)$ of all possible vertex pairs $<x'_{i},y'_{j}>$ and $<y'_{i},x'_{j}>$ is

\begin{equation}\label{Section-3-1-3}
\mathcal{S'}(3)=2\mathcal{S'}(1)-|\mathcal{\mathcal{V}}|(|\mathcal{V}|-1)
\end{equation}

Taken together cases 1.1-1.3, Eqs.(\ref{Section-3-1-1})-(\ref{Section-3-1-3}) produces an exact solution for $\mathcal{S'}$ parallel to that of Eq.(\ref{Section-3-1-0}) after some simple arithmetics. This completes our proof.

\emph{Corollary } After $t$ time steps, geodesic distance $\mathcal{S'}_{t}$ of first-order subdivision tree $\mathcal{T'}_{t}(\mathcal{V'}_{t},\mathcal{E'}_{t})$ will follow

\begin{equation}\label{Section-3-1-4}
\mathcal{S'}_{t}=8^{t}\mathcal{S}-\frac{1}{3}(2^{3t}-2^{t})(|\mathcal{V}|-1)+(2^{2t-1}-2^{3t-1})(|\mathcal{V}|-1)^{2}
\end{equation}

\textbf{Theorem 2} Given an arbitrary tree $\mathcal{T}(\mathcal{V},\mathcal{E})$, the exact solution for geodesic distance $\mathcal{S^{\odot}}$ of its $m$-vertex-operation tree $\mathcal{T^{\odot}}(\mathcal{V^{\odot}},\mathcal{E^{\odot}})$ is

\begin{equation}\label{Section-3-2-0}
\mathcal{S^{\odot}}=(1+m)^{2}\mathcal{S}+m(m+1)|\mathcal{V}|^{2}-m|\mathcal{V}|
\end{equation}
where $\mathcal{S}$ is a provided expression to geodesic distance of tree $T(\mathcal{V},\mathcal{E})$ initially.

\emph{\textbf{Proof}} By Def.3, $m$-vertex-operation tree $T^{\odot}(\mathcal{V^{\odot}},\mathcal{E^{\odot}})$ also contains two kinds of vertices, vertex sets $\mathcal{V}$ and $\mathcal{V^{\odot}}-\mathcal{V}$. To prevent new symbols from better understanding our proof, we here still make use of the above symbols in the process of developing theorem 1, that is, $X'=\mathcal{V}$ and $Y'=\mathcal{V^{\odot}}-\mathcal{V}$. As before, we will consider the different contribution from several various cases to the computation of geodesic distance $\mathcal{S^{\odot}}$ as follows.

\emph{Case 2.1} For a given vertex pair $<x'_{i},x'_{j}>$ of set $X'$, there must be a bijection $f_{1}$ between set $X'$ and set $\mathcal{V}$ such that $f_{1}(<x'_{i},x'_{j}>)=<x_{i},x_{j}>$ here vertices $x_{i}$ and $x_{j}$ are in set $\mathcal{V}$. In fact, such a bijection $f_{1}$ is self-mapping. Meanwhile, an $m$-vertex-operation is just applied to each vertex of tree $T(\mathcal{V},\mathcal{E})$ and hence has no influence on changing the geodesic distances $\mathcal{S}$ of all possible pair of vertices $x_{i}$ and $x_{j}$. So one can obtain

\begin{equation}\label{Section-3-2-1}
\mathcal{S^{\odot}}(1)=\mathcal{S}.
\end{equation}

\emph{Case 2.2} This is slightly different from discussion in case 1.2 of proof of theorem 1. We would like to bipartition the calculation of $\mathcal{S^{\odot}}(2)$ such that a portion of contribution comes from geodesic distances $\mathcal{S^{\odot}}(21)$, which is equal to the sum of geodesic distances on each pair of new vertices within each star, and the other $\mathcal{S^{\odot}}(22)$ is equivalent to the sum of geodesic distances on vertex pairs whose members are chosen from differen stars respectively. Compared to calculation of $\mathcal{S^{\odot}}(22)$, $\mathcal{S^{\odot}}(21)$ can be easily obtained and satisfies the below equation

\begin{equation}\label{Section-3-2-2}
\mathcal{S^{\odot}}(21)=|\mathcal{V}|\times \left(2\frac{m(m-1)}{2}\right).
\end{equation}

To calculate $\mathcal{S^{\odot}}(22)$, we need to construct a mapping $f_{2}$ between sets $Y'$ and $X'$. For a pair of vertices $y'_{i}$ connected to vertex $x'_{i}$ and $y'_{j}$ to $x'_{j}$, one must find out a unique path $P_{y'_{i}y'_{j}}$ which is made up edge set $y'_{i}x'_{i}$,..., $x'_{u}x'_{v}$,...,$x'_{j}y'_{j}$ and so there is a mapping $f^{*}$ between vertex pairs $<y'_{i},y'_{j}>$ and $<x'_{i},x'_{j}>$. Obviously, such mapping $f^{*}$ is a surjection because of vertex $x'_{i}$ adjacent to $m$ new vertices in view of $m$-vertex-operation. Therefore, we set mapping $f_{2}=f_{1}\circ f^{*}$ which is in fact an $m^{2}$-regular surjection and further have the following equation

\begin{equation}\label{Section-3-2-3}
\mathcal{S^{\odot}}(22)=m^{2}\mathcal{S^{\odot}}(1)+\sum_{i=1}^{|\mathcal{V}|-1}2m^{2}(|\mathcal{V}|-i).
\end{equation}

\emph{Case 2.3} For arbitrary vertex pair $<y'_{i},x'_{j}>$ where both subscripts may be the same, we can propose an $m$-regular surjection $f_{3}=f_{1}\circ f^{**}$ in which mapping $f^{**}$ maps vertex pair $<y'_{i},x'_{j}>$ to pair $<x'_{i},x'_{j}>$. In particular, mapping $f^{**}$ will map a vertex pair to a single vertex when $i=j$. Therefore, geodesic distance $\mathcal{S^{\odot}}(3)$ between arbitrary vertex $y'_{i}$ in set $Y'$ and each vertex $x'_{i}$ in set $X'$ should obey

\begin{equation}\label{Section-3-2-4}
\mathcal{S^{\odot}}(3)=m|\mathcal{V}|^{2}+2m\mathcal{S}.
\end{equation}

Combining Eqs.(\ref{Section-3-2-1})-(\ref{Section-3-2-4}) yields a precise solution for geodesic distance $\mathcal{S^{\odot}}$ equal to that of Eq.(\ref{Section-3-2-0}). This is complete.

From the appearance point of view, the ($1+m$)-star-fractal-operation on graph $\mathcal{G}(\mathcal{V},\mathcal{E})$ might also be achieved by both first-order subdivision and $m$-vertex-operation. Such an explanation is to first apply first-order subdivision to each edge of set $\mathcal{E}$ and then to manipulate $m$-vertex-operation only to those new vertices added by the first-order subdivision. Based on this close connection among them, the proof of theorem 3 can be concisely developed and completed. Now, let us divert our attention to turn out theorem 3.

\textbf{Theorem 3} Given arbitrary tree $\mathcal{T}(\mathcal{V},\mathcal{E})$, the exact solution for geodesic distance $\mathcal{S^{\star}}$ of its ($1+m$)-star-fractal-operation tree $\mathcal{T^{\star}}(\mathcal{V^{\star}},\mathcal{E^{\star}})$ is
\begin{equation}\label{Section-3-3-0}
\mathcal{S^{\star}}=2(m+2)^{2}\mathcal{S}-(m+2)(|\mathcal{V}|-1)(m+|\mathcal{V}|)
\end{equation}
here $\mathcal{S}$ is an already given expression to geodesic distance of tree $\mathcal{T}(\mathcal{V},\mathcal{E})$.

\emph{\textbf{Proof}} While there also are two different classes of vertices, vertex sets $X'=\mathcal{V}$ and $Y'=\mathcal{V^{\star}}-\mathcal{V}$, and three kinds of contributions to computation of geodesic distance $\mathcal{S^{\star}}$, it appears to be a little hard to capture closed-form of $\mathcal{S^{\star}}$ than the preceding two as we will show below. By analogy with the developments of the both theorems above, we will again make use of classification method which is more fine-grained than the foregoing both only because of the difference of topological structure among them.

\emph{Case 3.1} For a given vertex pair $<x'_{i},x'_{j}>$ of set $X'$, analogously, we can generate a self-mapping $f_{1}$ between set $X'$ and set $\mathcal{V}$ such that $f_{1}(<x'_{i},x'_{j}>)=<x_{i},x_{j}>$ here vertices $x_{i}$ and $x_{j}$ are in set $\mathcal{V}$. Therefore, geodesic distance $\mathcal{S^{\star}}(1)$ on such type of vertex pairs complies to

\begin{equation}\label{Section-3-3-1}
\mathcal{S^{\star}}(1)=2\mathcal{S}.
\end{equation}

\emph{Case 3.2} There will, by Def.2, be $|\mathcal{V}|-1$ new stars added into original tree $\mathcal{T}(\mathcal{V},\mathcal{E})$ in terms of ($1+m$)-star-fractal operation. Here just counts geodesic distances on arbitrary pair of leaf vertices within each of these stars but not between two distinct stars. This is in spirit similar to Eq.(\ref{Section-3-2-2}) and hence geodesic distance $\mathcal{S^{\star}}(2)$ on all possible leaf vertex pairs of this type is

\begin{equation}\label{Section-3-3-2}
\mathcal{S^{\star}}(2)=(|\mathcal{V}|-1)\left(2\frac{m(m-1)}{2}\right).
\end{equation}

\emph{Case 3.3} The majority of new vertices added by ($1+m$)-star-fractal-operation are leaf vertices, which constitute a set $Y'_{1}$ and the other portion of which, $Y'_{2}=Y'-Y'_{1}$, contains all vertices inserted on edges of tree $\mathcal{T}(\mathcal{V},\mathcal{E})$ as a whole. With the help of discussion mentioned in case 1.2, we still select bijection $f_{2}$, i.e., $f_{2}=f_{1}\circ f^{*}$, as our candidate mapping each element pair of set $Y'_{2}$ to a unique vertex pair of set $X'$. Under such a bijection, geodesic distance $\mathcal{S^{\star}}(3)$ on all possible vertex pairs of set $Y'_{2}$ should have the same outline as Eq.(\ref{Section-3-1-2})

\begin{equation}\label{Section-3-3-3}
\mathcal{S^{\star}}(3)=\mathcal{S^{\star}}(1)-2\frac{|\mathcal{V}|(|\mathcal{V}|-1)}{2}.
\end{equation}

\emph{Case 3.4} As shown in case 1.3, there have to be a surjection $f_{3}=f_{1}\circ f^{**}$ between vertex pairs $<x'_{i},y'_{j}>$, where $x'_{i}\in X'$, $y'_{j}\in Y'_{2}$ and subscripts $i$ and $j$ can be equal, and vertex pairs $<x'_{i},x'_{j}>$ in which $x'_{i}$ and $x'_{j}$ belong to set $X'$. In addition, there is in practice the other surjection between two such vertex pairs as will be stated here. For a given vertex pair $<x'_{i},y'_{j}>$, one can construct a mapping $f^{3*}$ mapping this pair to a unique vertex pair $<y'_{i},y'_{j}>$ of set $Y'_{2}$ by detecting an edge $x'_{i}y'_{i}$. As before, this mapping $f^{3*}$ is also a surjection according to another pre-image $<y'_{i},x'_{j}>$ following $(f^{3*})^{-1}(<y'_{i},y'_{j}>)=<y'_{i},x'_{j}>$. Therefore we adopt surjection $f^{3*}$ and bijection $f_{2}$ to generate our anticipated surjection $f_{3}$, i.e., $f_{3}=f_{2}\circ f^{3*}$, between vertex pair $<x'_{i},y'_{j}>$ and pair $<x'_{i},x'_{j}>$ and the closed-form formula for geodesic distance $\mathcal{S^{\star}}(4)$ on all possible vertex pairs $<x'_{i},y'_{j}>$ is

\begin{equation}\label{Section-3-3-4}
\mathcal{S^{\star}}(4)=2\mathcal{S^{\star}}(3)+|\mathcal{V}|(|\mathcal{V}|-1).
\end{equation}

\emph{Case 3.5} We here discuss geodesic distance $\mathcal{S^{\star}}(5)$ on all possible leaf vertex pairs $<y'_{i},y'_{j}>$ whose two vertices are from different stars not within a star. To accurately accomplish this, we have to produce a novel mapping $f^{4*}$, which is in essence an $m^{2}$-regular surjection, mapping vertex pairs $<y'_{i},y'_{j}>$ to vertex pairs $<y'_{io},y'_{jo}>$ here both vertices $y'_{io}$ and $y'_{jo}$ are, respectively, the central vertex of stars to which vertices $y'_{i}$ and $y'_{j}$ belong. Using that bijection $f_{2}$ introduced in case 3.3, we can propose a reasonable surjection $f_{4}=f_{2}\circ f^{4*}$ at once and write

\begin{equation}\label{Section-3-3-5}
\mathcal{S^{\star}}(5)=m^{2}\mathcal{S^{\star}}(3)+\sum_{i=1}^{|\mathcal{V}|-2}2m^{2}(|\mathcal{V}|-1-i).
\end{equation}

\emph{Case 3.6} Now let us focus on computation of geodesic distance $\mathcal{S^{\star}}(6)$ on all possible vertex pairs $<x'_{i},y'_{j}>$ where $x'_{i}\in X'$, $y'_{j}\in Y'_{1}$ and subscripts $i$ and $j$ can be the same. Considering such a vertex pair $<x'_{i},y'_{j}>$, one can first find out a surjection $f^{5*}$ such that $f^{5*}(<x'_{i},y'_{j}>)=<x'_{i},y'_{jo}>$ where vertex $y'_{jo}$ is the central vertex of star to which vertex $y'_{j}$ belongs. And then, taking well constructed surjection $f_{3}$ in case 3.4, we can establish an acceptable surjection $f_{5}=f_{3}\circ f^{5*}$ between vertex pairs $<x'_{i},y'_{j}>$ and $<x'_{i},y'_{jo}>$ and thus the solution for geodesic distance $\mathcal{S^{\star}}(6)$ obeys

\begin{equation}\label{Section-3-3-6}
\mathcal{S^{\star}}(6)=m\mathcal{S^{\star}}(4)+m|\mathcal{V}|(|\mathcal{V}|-1).
\end{equation}

\emph{Case 3.7} So far, we have finished the entire computations of geodesic distance on vertex pairs whose one vertex is in set $X'$ and the other in set $Y'$. The left is to measure geodesic distance $\mathcal{S^{\star}}(7)$ on vertex pairs in which one vertex is chosen from set $Y'_{1}$ and another one from set $Y'_{2}$. With the above notations, such a vertex pair can be regarded as $<y'_{i},y'_{jo}>$ in which it is likely that subscript $i$ is equivalent to $j$. For all but $i=j$ vertex pairs $<y'_{i},y'_{jo}>$, it is natural to present an $m$-regular surjection $f^{6*}$ pointing vertex pair $<y'_{i},y'_{jo}>$ to $<y'_{io},y'_{jo}>$ which together with bijection $f_{2}$ in case 3.3 both are compounded into our expectant surjection $f_{6}=f_{2}\circ f^{6*}$ that is able to do what we want to do. Therefore, an expression of geodesic distance $\mathcal{S^{\star}}(7)$ can be written as

\begin{equation}\label{Section-3-3-7}
\mathcal{S^{\star}}(7)=2m\mathcal{S^{\star}}(3)+m(|\mathcal{V}|-1)^{2}.
\end{equation}

Substituting Eqs.(\ref{Section-3-3-1})-(\ref{Section-3-3-7}) into this expression $\sum_{i=1}^{7}\mathcal{S^{\star}}(i)$ and adopting some simple arithmetics together produces the same result as shown in Eq.(\ref{Section-3-3-0}) which is complete.

\section{Applications}

In this section, we will introduce two families of well studied treelike models, $\mathcal{T}^{\odot}(t,m)$ and $\mathcal{T}^{\star}(t,m+1)$, and then apply our methods for calculating geodesic distances on the two models. We are not the first to study geodesic distance on both models $\mathcal{T}^{\odot}(t,m)$ and $\mathcal{T}^{\star}(t,m+1)$. There are some published papers reporting such problems in the last few decades due to the two models themselves have some other interesting topological structure properties. The model $\mathcal{T}^{\odot}(t,m)$, for instance, obeys exponential degree distribution in form and on the other hand model $\mathcal{T}^{\star}(t,m+1)$ exhibits in nature fractal feature. At the same time, the both have self-similar character. It is indeed on this basis of self-similarity that many researchers had captured exact solutions for geodesic distance on models $\mathcal{T}^{\odot}(t,m)$ and $\mathcal{T}^{\star}(t,m+1)$. Some of which all have a common assumption that the original graph $\mathcal{G}(\mathcal{V},\mathcal{E})$ is a single edge connecting two vertices. Apparently, these methods based on such an assumption have most likely to become tough and intractable when the seed graph $\mathcal{G}(\mathcal{V},\mathcal{E})$ is arbitrary one we are interested in. To address this issues, many take useful advantage of algebraic methods from spectral graph theory, such as eigenvalue and laplacian matrix. As known, facing with some larger graphs with a great number of vertices designed as seeds, methods built by matrix will also become cumbersome by respect with limitation of space and memory on computers themselves. Our methods may, by contrast, play a better role on deriving accurate solutions for some treelike models of specifical types as we will discuss shortly.

As the first example of applications of our methods to some treelike models, model $\mathcal{T}^{\odot}(t,m)$ cam be certainly generated using the next algorithm.

\subsection{Model $\mathcal{T}^{\odot}(t,m)$}

Here we also choose a single edge as a seed to build up the whole treelike model $\mathcal{T}^{\odot}(t,m)$. One of the most important reasons for this is to make a comparison with some known methods to highlight conciseness and convenience of our methods, at least on this model.

\textbf{\emph{Algorithm 1}}

At $t=0$, the seminal graph $\mathcal{T}$ is indeed an edge incident to two vertices, shown in Fig.2(a), which is also viewed as $\mathcal{T}^{\odot}(0,m)$ for convenience.

At $t=1$, the next model $\mathcal{T}^{\odot}(1,m)$ can be obtained from model $\mathcal{T}^{\odot}(0,m)$ by manipulating $m$-vertex-operation on each vertex as illustrated in Fig.2(b) where parameter $m=2$.

At $t=2$, the second model $\mathcal{T}^{\odot}(2,m)$ can be generated in a similar manner to that mentioned at the preceding time step, shown in Fig.2(c) where parameter $m=2$ as well.

For $t\geq3$, algorithm 1 can well run to output a model $\mathcal{T}^{\odot}(t,m)$ what we are interested in. Before beginning by our discussion, let we first report one of the most prominent topological structure properties of model $\mathcal{T}^{\odot}(t,m)$.

\begin{figure}
\centering
  \includegraphics[height=6cm]{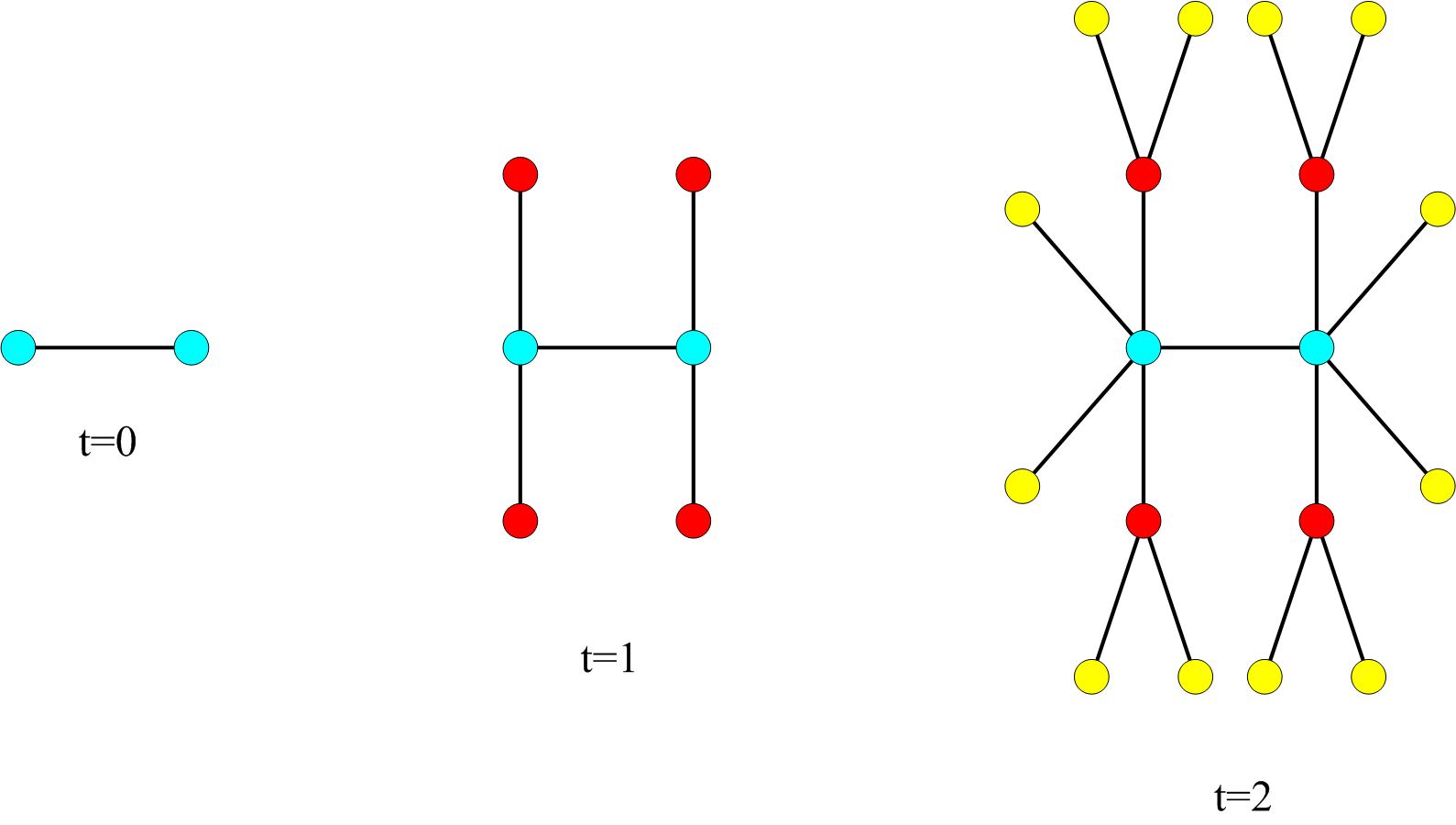}\\
{\small Fig.2. The diagram of the first three of the graph-model $\mathcal{T}^{\odot}(t,m)$ where $m=2$ .       }
\end{figure}

It can be seen that model $\mathcal{T}^{\odot}(t,m)$ has $t+1$ different types of vertices according to vertex degree. Expect for that two vertices attached to original edge, the total number of vertices added at time step $t_{i}$ into model $\mathcal{T}^{\odot}(t,m)$ is equal to $2m(m+1)^{t_{i}-1}$ and each of them has the same degree $k_{t_{i}}=2(t-t_{i}+1)$. By the definition of cumulative degree distribution, model $\mathcal{T}^{\odot}(t,m)$ obeys

\begin{equation}\label{Section-4-1-0}
P_{cum}(k\geq k_{t_{i}})\sim \emph{exp}[-(m+1)k_{t_{i}}].
\end{equation}
Here states an exponential degree distribution with parameter $-(m+1)$.

Our reasons for choosing this model $\mathcal{T}^{\odot}(t,m)$ with exponential degree distribution are twofold. First of all, such a degree distribution can be prevailing in a great deal of small-world network models
including the well known WS model \cite{Watts-1998}. The other is that based on this degree distribution one can build up many other network models of interest, for instance, power-law networks according to some quantity distribution such as weight. To see why this is so, we here introduce an illustrated example as follows. Consider a network model with exponential degree distribution, that is, the degree value $k$ of vertex follows an exponential distribution $P(k)\sim$ exp$(\alpha k)$ in which parameter $\alpha$ indicates whether this distribution is normalizable or not. Now, we can be able to assume, in some case, that the weight value $w_{k}$ of degree $k$ vertex is exponentially related to degree value $k$, i.e., $w_{k}\sim$ exp$(\beta k)$ where symbol $\beta$ represents another parameter. And then the probability distribution of weight $w_{k}$ obeys the following expression
\begin{equation}\label{Section-4-1-1}
P(w_{k})=P(k)\frac{dk}{dw_{k}}\sim\frac{w_{k}^{-1+\frac{\alpha}{\beta}}}{\beta}.
\end{equation}
Obviously, this is a power law with exponent $\gamma=1-\frac{\alpha}{\beta}$. Hence, it is of prominent interest to study graphs having exponential degree distribution associated with some quantities in terms of both realistic and theoretical senses.

Now let us put insight to computation of geodesic distance $\mathcal{S^{\odot}}_{t}$ on model $\mathcal{T}^{\odot}(t,m)$. With the description of algorithm 1, model $\mathcal{T}^{\odot}(t,m)$ can be generated from model $\mathcal{T}^{\odot}(t-1,m)$ using $m$-vertex-operation, model $\mathcal{T}^{\odot}(t-1,m)$ from model $\mathcal{T}^{\odot}(t-2,m)$, and so forth. An iterative equation may be established based on Eq.(\ref{Section-3-2-0}) in the following way

\begin{equation}\label{Section-4-1-2}
\mathcal{S^{\odot}}_{t}=(1+m)^{2}\mathcal{S^{\odot}}_{t-1}+m(m+1)|\mathcal{V^{\odot}}_{t-1}|^{2}-m|\mathcal{V^{\odot}}_{t-1}|.
\end{equation}

Adopting some simple arithmetics, we can reorganized Eq.(\ref{Section-4-1-2}) as

\begin{equation}\label{Section-4-1-3}
\mathcal{S^{\odot}}_{t}=(m+1)^{2t}\mathcal{S^{\odot}}_{0}+m(m+1)\sum_{i=0}^{t-1}(m+1)^{2i}|\mathcal{V^{\odot}}_{t-1-i}|^{2}-m\sum_{i=0}^{t-1}(m+1)^{2i}|\mathcal{V^{\odot}}_{t-1-i}|.
\end{equation}

With some additional conditions $\mathcal{S^{\odot}}_{0}=1$ and $|\mathcal{V^{\odot}}_{t}|=2(m+1)^{t}$, an exact solution for geodesic distance $\mathcal{S^{\odot}}_{t}$ is

\begin{equation}\label{Section-4-1-4}
\mathcal{S^{\odot}}_{t}=(m+1)^{t-1}[2+(4mt+m-1)(m+1)^{t}].
\end{equation}

Meanwhile, we can capture a closed-form expression of average geodesic distance $\langle\mathcal{S^{\odot}}_{t}\rangle$ on model  $\mathcal{T}^{\odot}(t,m)$ as follows

\begin{equation}\label{Section-4-1-5}
\langle\mathcal{S^{\odot}}_{t}\rangle=\frac{\mathcal{S^{\odot}}_{t}}{|\mathcal{V^{\odot}}_{t}|(|\mathcal{V^{\odot}}_{t}|-1)/2}\sim 2t
\end{equation}
in the limit of large graph size, shown in Fig.3.

\begin{figure}
\centering
  \includegraphics[height=6cm]{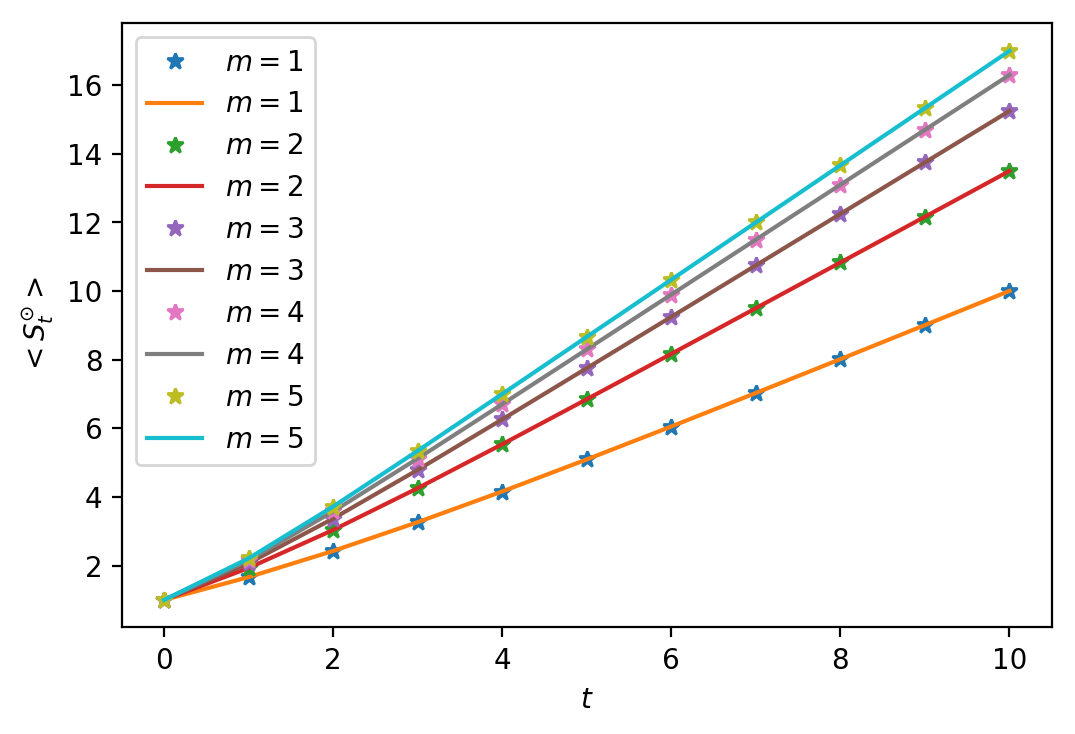}\\
{\small Fig.3. The diagram of mean geodesic distance $\langle\mathcal{S^{\odot}}_{t}\rangle$ on graph-model $\mathcal{T}^{\odot}(t,m)$. The solid lines represent analytical values serving as guides to the eye and the stars indicate computer simulations.         }
\end{figure}

In addition, the diameter $\mathcal{D}^{\odot}(t,m)$ of model $\mathcal{T}^{\odot}(t,m)$ is $2t+1$ and the logarithm value of vertex number $|\mathcal{V^{\odot}}_{t}|$ is asymptotically equal to $t\ln(m+1)$. Plugging the two values into Eq.(\ref{Section-4-1-5}) yields

\begin{equation}\label{Section-4-1-6}
\langle\mathcal{S^{\odot}}_{t}\rangle\sim \mathcal{D}^{\odot}(t,m)\sim \ln |\mathcal{V^{\odot}}_{t}|
\end{equation}
in the large graph size limit. From the complexity point of view, diameter and logarithm value of vertex number both perform better than average geodesic distance when we attempt to measure how long a pair of information throughout diffuses on the entire model $\mathcal{T}^{\odot}(t,m)$. This indicates that model $\mathcal{T}^{\odot}(t,m)$ has smaller diameter and lower average geodesic distance comparison with vertex number. Perhaps this is why more and more researchers would like to select diameter or logarithm value of vertex number rather than average geodesic distance to determine whether a complex network under consideration is small-world or not.

Last but not least, the mean first-passage time $\mathcal{R}^{\odot}(t,m)$ for a random walker on model $\mathcal{T}^{\odot}(t,m)$ can be expressed in the below equation

\begin{equation}\label{Section-4-1-7}
\mathcal{R}^{\odot}(t,m)=\frac{2\mathcal{S^{\odot}}_{t}}{|\mathcal{V^{\odot}}_{t}|}=\frac{2}{m+1}+(4mt+m-1)(m+1)^{t-1}
\end{equation}
which is completely parallel to that of Ref \cite{Zhang-2010} but our methods are quite light and concise compared to the latter.

We want to stress again that in the procedure of developing model $\mathcal{T}^{\odot}(t,m)$ the seed may be any graph of importance interest. When the seed is a relatively large one, some published methods, such as that addressed in \cite{Zhang-2010}, can not be chosen as an adequate candidate for calculating geodesic distance precisely. However, our methods seem to be a potential scheme for addressing such problems, at least from the aspect of complexity. We look forward to seeing some applications of our methods to other models in the next future.

Another potential application of our methods is to derive a precise solution for geodesic distance $\mathcal{S^{\star}}(t,m+1)$ on a best studied treelike model $\mathcal{T}^{\star}(t,m+1)$ because this model possesses some interesting topological structure properties, for instance, fractal character, and is in nature solvable. Below is a construction of model $\mathcal{T}^{\star}(t,m+1)$ by virtue of algorithm 2.

\subsection{Model $\mathcal{T}^{\star}(t,m+1)$}

As before, we still choose a single edge as a seed to construct the whole treelike model $\mathcal{T}^{\star}(t,m+1)$. In principle, the seed can be an arbitrary tree as stated previously.

\textbf{\emph{Algorithm 2}}

At $t=0$, the seminal graph $\mathcal{T}$ is indeed an edge incident to two vertices, shown in Fig.4(a), which is also thought of as $\mathcal{T}^{\star}(0,m+1)$ for convenience.

At $t=1$, the newborn model $\mathcal{T}^{\star}(1,m+1)$ can be generated from model $\mathcal{T}^{\star}(0,m+1)$ by executing ($1,m$)-star-fractal operation on each edge as plotted in Fig.4(b) where parameter $m=2$.

At $t=2$, the second model $\mathcal{T}^{\star}(2,m+1)$ can be produced in a same manner as that of adopted at the foregoing time step, illustrated in Fig.4(c) where parameter $m=2$ also.

For $t\geq3$, one can terminate algorithm 2 until he (she) obtains a desired model $\mathcal{T}^{\star}(t,m+1)$ in which parameter $t$ is designed in advance. Before beginning by our discussion, let we first enumerate the total number $|\mathcal{V^{\star}}_{t}|$ of vertices of model $\mathcal{T}^{\star}(t,m+1)$ in a recursive way

\begin{equation}\label{Section-4-2-0}
|\mathcal{V^{\star}}_{t}|=|\mathcal{V^{\star}}_{t-1}|+(m+1)(|\mathcal{V^{\star}}_{t-1}|-1).
\end{equation}

Solving for $|\mathcal{V^{\star}}_{t}|$ with initial condition $|\mathcal{V^{\star}}_{0}|=2$ produces $|\mathcal{V^{\star}}_{t}|=(m+2)^{t}(|\mathcal{V}|-1)+1$ where $|\mathcal{V}|$ is vertex number of seminal graph $\mathcal{T}$ and in fact equal to $2$ and hence $|\mathcal{V^{\star}}_{t}|=(m+2)^{t}+1$ finally.

\begin{figure}
\centering
  \includegraphics[height=6cm]{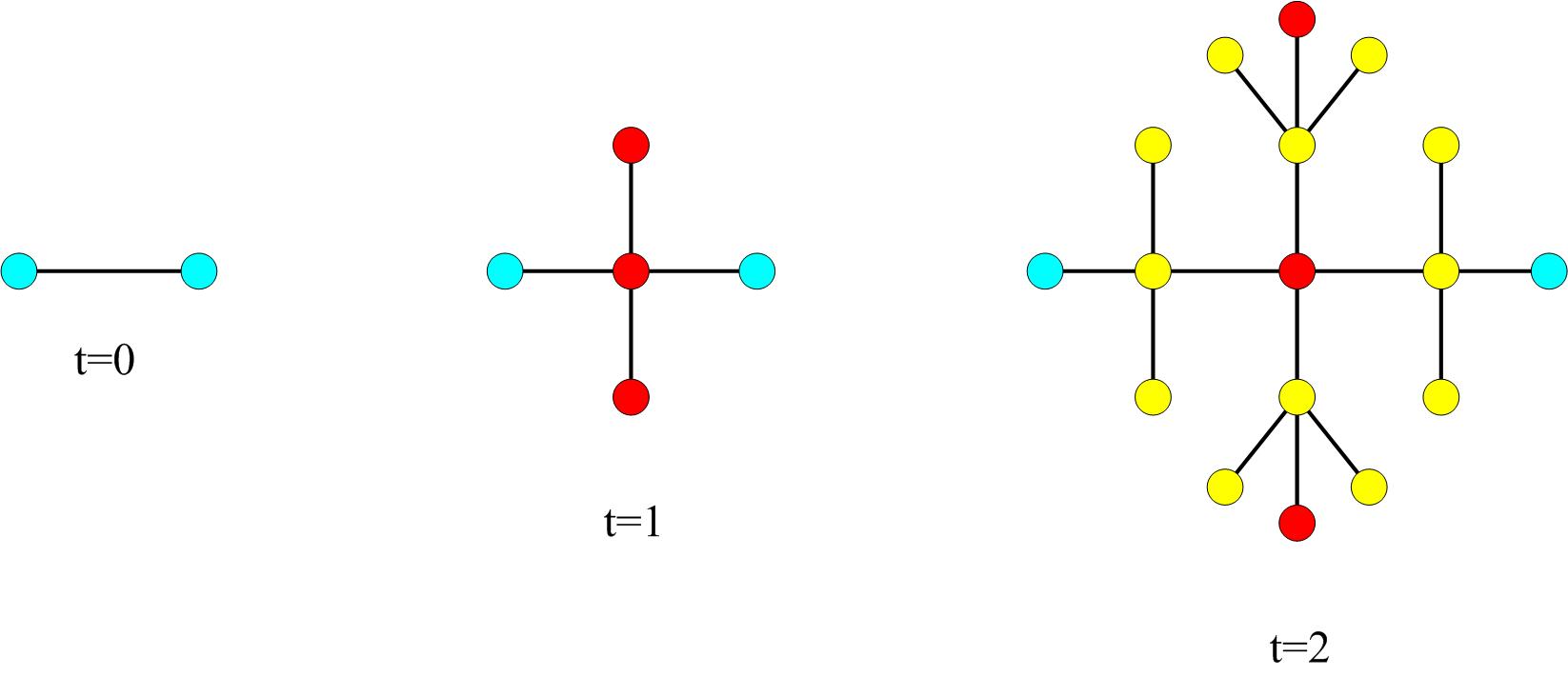}\\
{\small Fig.4. The diagram of the first three of the graph-model $\mathcal{T}^{\star}(t,m+1)$ where $m=2$.       }
\end{figure}

From now on, we are going to calculate geodesic distance $\mathcal{S^{\star}}_{t}$ on model $\mathcal{T}^{\star}(t,m+1)$ in a fashion similar to that employed to build up Eq.(\ref{Section-4-1-2}). It is clear to the eye that model $\mathcal{T}^{\star}(t,m+1)$ is certainly established by applying ($1,m$)-star-fractal operation to each edge of the ancestor model $\mathcal{T}^{\star}(t-1,m+1)$, model $\mathcal{T}^{\star}(t-1,m+1)$ established by its ancestor $\mathcal{T}^{\star}(t-2,m+1)$, and so on. This allows us to capture an accurate solution for geodesic distance $\mathcal{S^{\star}}_{t}$ in an iterative way. With the enlightenment from the proof of theorem 3, we arrive at

\begin{equation}\label{Section-4-2-1}
\mathcal{S^{\star}}_{t}=2(m+2)^{2}\mathcal{S^{\star}}_{t-1}-(m+2)(|\mathcal{V^{\star}}_{t-1}|-1)(m+|\mathcal{V^{\star}}_{t-1}|).
\end{equation}

Substituting Eq.(\ref{Section-4-2-0}) into Eq.(\ref{Section-4-2-1}) outputs

\begin{equation}\label{Section-4-2-2}
\begin{aligned}\mathcal{S^{\star}}_{t}&=[2(m+2)^{2}]^{t}\mathcal{S^{\star}}_{0}-(m+2)\sum_{i=0}^{t-1}[2(m+2)^{2}]^{i}|\mathcal{V^{\star}}_{t-1-i}|^{2}+m(m+2)\sum_{i=0}^{t-1}[2(m+2)^{2}]^{i}\\
&-(m-1)(m+2)\sum_{i=0}^{t-1}[2(m+2)^{2}]^{i}|\mathcal{V^{\star}}_{t-1-i}|.
\end{aligned}
\end{equation}

With the initial conditions $\mathcal{S^{\star}}_{0}=1$ and $|\mathcal{V^{\star}}_{0}|=2$, an exact solution for geodesic distance $\mathcal{S^{\star}}_{t}$ follows

\begin{equation}\label{Section-4-2-3}
\mathcal{S^{\star}}_{t}=[2(m+2)^{2}]^{t}-(2^{t}-1)(m+2)^{2t-1}-\frac{(m+1)(m+2)^{t}([2(m+2)]^{t}-1)}{2(m+2)-1}.
\end{equation}

At the same time, it is easy to obtain a closed-form expression of average geodesic distance $\langle\mathcal{S^{\star}}_{t}\rangle$ on model $\mathcal{T}^{\star}(t,m+1)$

\begin{equation}\label{Section-4-2-4}
\langle\mathcal{S^{\star}}_{t}\rangle=\frac{\mathcal{S^{\star}}_{t}}{|\mathcal{V^{\star}}_{t}|(|\mathcal{V^{\star}}_{t}|-1)/2}\sim 2^{t}
\end{equation}
in the limit of large graph size, plotted in Fig.5.

\begin{figure}
\centering
  \includegraphics[height=6cm]{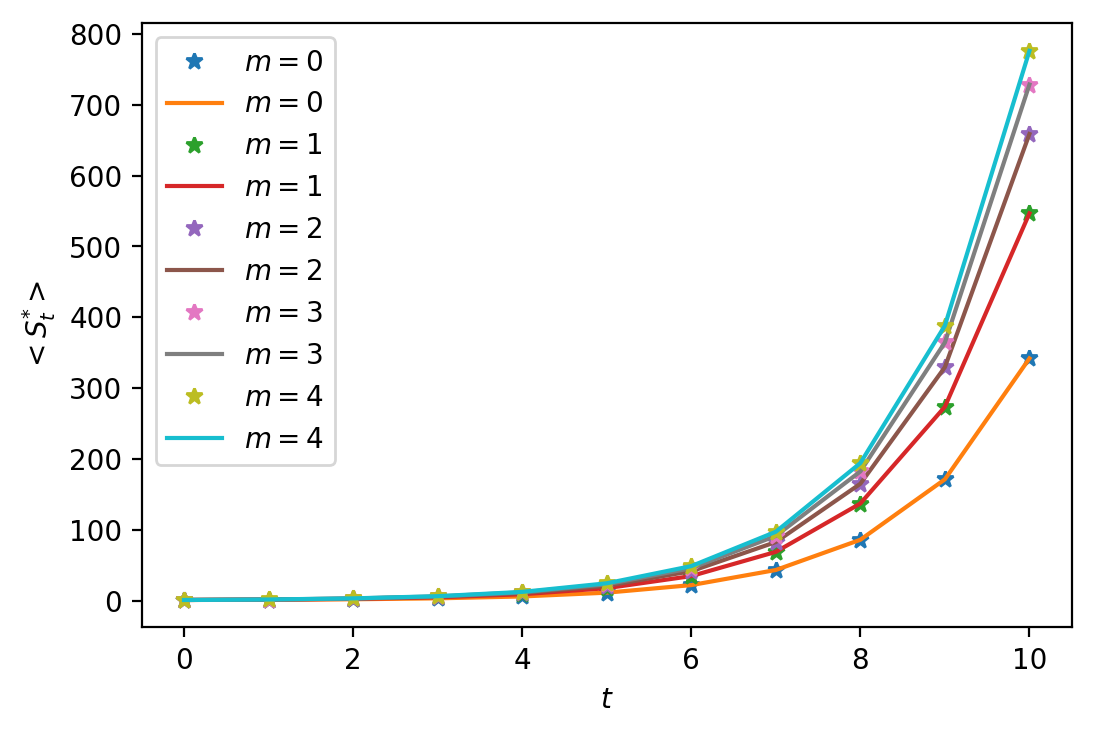}\\
{\small Fig.5. The diagram of mean geodesic distance $\langle\mathcal{S^{\star}}_{t}\rangle$ on graph-model $\mathcal{T}^{\star}(t,m+1)$. The solid lines represent analytical values serving as guides to the eye and the stars indicate computer simulations.         }
\end{figure}

As discussed above, there is a relationship among average geodesic distance $\langle\mathcal{S^{\star}}_{t}\rangle$, diameter $\mathcal{D}^{\star}(t,m+1)$ and the logarithm value of vertex number $\ln|\mathcal{V^{\star}}_{t}|$ as below

\begin{equation}\label{Section-4-2-5}
\langle\mathcal{S^{\star}}_{t}\rangle\sim \mathcal{D}^{\star}(t,m+1)\gg \ln |\mathcal{V^{\star}}_{t}|
\end{equation}
in the large graph size limit. Different from Eq.(\ref{Section-4-1-6}),  Eq.(\ref{Section-4-2-5}) manifests a fact that model $\mathcal{T}^{\star}(t,m+1)$ is large-size. One of reasons for this is the considerable difference between $m$-vertex-operation and ($1,m$)-star-fractal operation which is in nature evident. In addition, model $\mathcal{T}^{\star}(t,m+1)$ is more homogeneous than model $\mathcal{T}^{\odot}(t,m+1)$ because all but leaf vertices have an identical smaller degree in the model as a whole.

By analogy with technique addressed in Eq.(\ref{Section-4-1-7}), the mean first-passage time $\mathcal{R}^{\star}(t,m+1)$ for a random walker on model $\mathcal{T}^{\star}(t,m+1)$ can be written as

\begin{equation}\label{Section-4-2-6}
\mathcal{R}^{\star}(t,m+1)=\frac{2\mathcal{S^{\star}}_{t}}{|\mathcal{V^{\star}}_{t}|}=2\times\frac{[2(m+2)^{2}]^{t}-(2^{t}-1)(m+2)^{2t-1}-\frac{(m+1)(m+2)^{t}([2(m+2)]^{t}-1)}{2(m+2)-1}}{(m+2)^{2}+1}
\end{equation}
which is completely equivalent to that of Ref \cite{Zhang-2011} but our methods are quite light and convenient when considering a larger tree as seed.

One special case of our model $\mathcal{T}^{\star}(t,m+1)$ has in fact taken much attention in some other fields \cite{O-2010}. This case is the so-called \textbf{\emph{T-graph}} which can be induced from model $\mathcal{T}^{\star}(t,m+1)$ by setting parameter $m=1$. We here review some prominent topological structure properties on T-graph briefly. The T-graph is a fractal with the fractal dimension $d_{f}=\frac{\ln3}{\ln2}$ and the random-walk dimension $d_{w}=\frac{\ln6}{\ln2}=1+d_{f}$. In the meantime, the spectral dimension of T-graph is $\widetilde{d}=\frac{2d_{f}}{d_{w}}=\frac{\ln 9}{\ln 6} <2$, suggesting which a random walk on it is persistent \cite{Erik-2005}. Similarly, our model $\mathcal{T}^{\star}(t,m+1)$ has the fractal dimension $d^{\star}_{f}=\frac{\ln(m+2)}{\ln2}$, the random-walk dimension $d^{\star}_{w}=\frac{\ln2(m+2)}{\ln2}=1+d_{f}$ and the spectral dimension $\widetilde{d}^{\star}=\frac{2d^{\star}_{f}}{d^{\star}_{w}}=\frac{\ln (m+2)^{2}}{\ln 2(m+2)} <2$. Therefore, a random walk on our model $\mathcal{T}^{\star}(t,m+1)$ is persistent as well.

With both Eq.(\ref{Section-4-1-7}) and Eq.(\ref{Section-4-2-6}), it is straightforward to say that the topological structure of underlying model has a significant influence on the mean first-passage time for a random walker on it. This indicates that in general the mean first-passage time for a random walker on one heterogeneous model is more smaller than that on a homogeneous one.

\section{Conclusion and problem}

In conclusion, we propose some simple and suitable methods, which are established by virtue of two useful mathematical tools, i.e., mapping and vertex cover, for computation of geodesic distance on many types of treelike models which are built by three kinds of operations, first-order subdivision, ($1,m$)-star-fractal operation and $m$-vertex-operation. While these models discussed in this paper have been well studied in the last decades and a number of topological structure properties of interest have been reported, including computation of geodesic distance, our methods for geodesic distance on these such models will become quite convenient compared with published ones when in general considering an arbitrary tree as a seed, in particular large-size one without self-similarity. Such advantages of our methods can be clear in the process of development of theorems.

To show the importance and some potential applications of our methods, we introduce two families of treelike models with respect to their own specific topological structures where the one obeys exponential degree distribution and the other has fractal feature. With methods addressed here, we can derive exact solutions for geodesic distance on the both models. In fact, these closed-form formulas are not fresh and already obtained by other authors in some published papers. Even so, this does not erase our contribution. On the one hand, these results captured here can serve as strong proofs to prove those published ones correct. Our methods, on the other hand, are able to be considered concise because of reducing some complicated calculations in comparison with many commonly used methods. One of the most important reasons for this is that we present these methods by taking into account structure features of the both treelike models throughout \cite{Fei-2018,Chris-2019}. Meantime, we have determined the formula of geodesic distance on $n$-order subdivision tree models in \cite{Fei-2019}.

We would like to stress that our work is only a tip of the iceberg and however the lights shed by our methods can be useful. We wish to witness some other applications of our methods in the days to come.

\textbf{Acknowledgment.} The authors would like to thank Bing Yao for useful conversations. This work was supported in part by the National Natural Science Foundation of China under grant No. 61662066.

\end{document}